\documentclass[12pt]{article}
\usepackage{amsfonts, amssymb, amsmath, mathrsfs, amsthm, latexsym, graphicx}
\usepackage{indentfirst}
\theoremstyle{plain}
\newtheorem{theorem}{Theorem}[section]

\newtheorem{corollary}[theorem]{Corollary}
\newtheorem{proposition}[theorem]{Proposition}

\theoremstyle{remark}

\theoremstyle{definition}

\numberwithin{equation}{section} 
\begin{document}
\title{Antimagic labelings of a complete graph}
\author{Dr.\ A.\ N.\ Bhavale}
\date{}
\maketitle
\begin{center}
hodmaths@moderncollegepune.edu.in \\ Department of Mathematics, \\ Modern College of Arts, Science and Commerce (Autonomous), \\ 
Shivajinagar, Pune 411005, M.S., India.
\end{center}
\begin{abstract}
In $1990$, Hartsfield and Ringel introduced antimagic graphs.  
Hartsfield and Ringel conjectured that every connected graph (and in particular, a tree) except $K_2$ is antimagic.
In $2010$, Hefetz et al.\ raised two questions: Is every orientation of any simple connected undirected graph antimagic? and
Given any undirected graph $G$, does there exist an orientation of $G$ which is antimagic? They call such an orientation an {\it antimagic orientation} of $G$.
Recently, Bhavale provided an edge labeling for a given graph on $n$ vertices without isolated vertices.
In this paper, using the labeling of Bhavale, we prove that a complete graph $K_n$ for $n \geq 3$ is super antimagic as well as totally antimagic total graph.
We also prove that there exists an antimagic orientation of $K_n$ for $n \geq 3$. 
\end{abstract}
MSC Classification $2020$: $05C78$.

\section{Introduction} \label{sec1}

According to Gallian \cite{JAG}, a {\it graph labeling} is an assignment of integers to the vertices or edges, or both, subject to certain conditions. 
A {\it{labeled graph}} $G=(V, E)$ is a finite series of graph vertices $V$ with a set of graph edges $E$ of $2$-subsets of $V$ (see Weisstein \cite{EW}). 
Formally, given a graph $G = (V, E)$, {\it a vertex labeling} is a function of $V$ to a set of labels; a graph with such a function defined is called {\it a vertex-labeled graph}. Likewise, {\it an edge labeling} is a function of $E$ to a set of labels. In this case, the graph is called {\it an edge-labeled graph}.
A {\it complete graph} on $n$ vertices, denoted by $K_n$, is a (simple) graph in which each vertex is connected to every other vertex by an edge. 

Graph labelings were first introduced in the mid $1960$s. The origin of graph labelings can be attributed to Rosa \cite{AR}. A general survey of graph labelings is found in Gallian \cite{JAG}. 
In the intervening years over $200$ graph labelings techniques have been studied in over $3000$ papers. There are many researchers from all over the world working on different kinds of graph labelings. In $1990$, Hartsfield and Ringel \cite{HR} introduced antimagic graphs, i.e., graphs having an antimagic labeling.  
Hartsfield and Ringel \cite{HR} conjectured that every tree as well as a connected graph except $K_2$ is antimagic.
Recently, Bhavale \cite{ANB} obtained an edge labeling for a given vertex labeling of a simple graph on $n$ unisolated vertices.
In this paper, using the labeling of Bhavale, we prove that a complete graph $K_n$ for $n \geq 3$ is super antimagic as well as totally antimagic total graph.
We also prove that there exists an antimagic orientation of $K_n$ for $n \geq 3$. 
For the other definitions, notation, and terminology see \cite{{HR},{DBW}}. 

\section{Antimagic labelings of graphs/directed graphs} \label{sec2}

Let $G=(V, E)$ be a graph with $|V| = n$ and $|E| = l$.
An {\it antimagic labeling} of a graph $G$ is a bijection $f : E \to \{1, 2, \ldots , l\}$ such that all the $n$ vertex sums are pairwise distinct, where a {\it vertex sum} is the sum of labels of all edges incident to that vertex.
Thus, if a vertex sum of a vertex $u$ is denoted by $S_u$ then $S_u = \displaystyle \sum_{e \in I(u)} f(e)$, where $I(u)$ is the set of all edges incident to $u$.
A graph is called {\it antimagic} if it has an antimagic labeling.
Hartsfield and Ringel \cite{HR} showed that paths $P_n$ ($n \geq 3$), cycles, wheels, and complete graphs $K_n$ ($n \geq 3$) are antimagic.

In $2010$, Hefetz et al.\ \cite{HMS} investigated antimagic labelings of directed graphs. An {\it antimagic labeling of a directed graph} $D$ with $n$ vertices and $l$ arcs (or directed edges) is a bijection from the set of arcs of $D$ to the integers $\{1, 2, \ldots , l\}$ such that all $n$ oriented vertex sums are pairwise distinct, where an {\it oriented vertex sum} is the sum of labels of all edges entering that vertex minus the sum of labels of all edges leaving it. 

We define an {\it antimagic labeling of a directed graph} $\overrightarrow{D}$ with vertices $\{v_1, v_2, \ldots, v_n\}$ and $l$ arcs (or directed edges) as a bijection from the set of arcs of $\overrightarrow{D}$ to the integers $\{1, 2, \ldots , l\}$ such that all $n$ vertex sums are pairwise distinct, where a {\it vertex sum} is the sum of labels of all edges entering as well as leaving that vertex. 
Let $S_{v_i}^- = \displaystyle \sum_{\overrightarrow{e} \in N^-(v_i)} f(\overrightarrow{e})$, where 
$N^-(v_i) = \{\overrightarrow{(v_p,v_i)} \in E(\overrightarrow{D}) | p < i \}$. Similarly, 
let $S_{v_i}^+ = \displaystyle \sum_{\overrightarrow{e} \in N^+(v_i)} f(\overrightarrow{e})$, where 
$N^+(v_i) = \{\overrightarrow{(v_i,v_q)} \in E(\overrightarrow{D}) | i < q \}$.
Then the vertex sum of vertex $v_i$ is given by $S_{v_i} = S_{v_i}^- + S_{v_i}^+$, and the oriented vertex sum of vertex $v_i$ is given by $S^o_{v_i} = S_{v_i}^- - S_{v_i}^+$.
A digraph is called {\it antimagic} if it has an antimagic labeling.
Let $S = \{(i, j) | 1 \leq i < j \leq n\}$.
Let $J_N = \{k | 1 \leq k \leq N = \binom{n}{2}\}$ where $n \geq 2$. 
The following result is due to Bhavale \cite{ANB}.

\begin{theorem} \cite{ANB} \label{t1}
Let $n \geq 2$. Define $F : S \to J_N$ by $F(i,j) = (i-1)n - \binom{i}{2} + j - i, \forall (i,j) \in S$. Then $F$ is bijective, and hence $|S| = N = \binom{n}{2}$.
\end{theorem}
Let $V(K_n) = \{v_1, v_2, \ldots , v_n\}$.
Let $\overrightarrow{K_n}$ be the directed graph associated with $K_n$ such that whenever $1 \leq i < j \leq n$ there is a directed edge from $v_i$ to $v_j$, (that is, $\overrightarrow{(v_i, v_j)} \in E(\overrightarrow{K_n})$).
Then the following two results follows from Theorem \ref{t1}. 

\begin{corollary} \cite{ANB} \label{c1}
Consider a directed graph $\overrightarrow{K_n}$ with vertex labeling $\{1, 2, \ldots , n\}$. Then for $i < j$, a directed edge $\overrightarrow{(i,j)} \in E(\overrightarrow{K_n})$ can be labeled (uniquely) as $k = (i-1)n - \binom{i}{2} + j - i \in J_N$. Moreover, for any $k \in J_N$ there exists unique directed edge $\overrightarrow{(i,j)} \in E(\overrightarrow{K_n})$.
\end{corollary}

\begin{corollary} \cite{ANB} \label{t2}
Suppose $\overrightarrow{G}$ is a directed subgraph (on $n$ unisolated vertices) of $\overrightarrow{K}_n$ with vertex labeling $\{v_1, v_2, \ldots, v_n\}$. Then for $i<j$, if $\overrightarrow{(v_i,v_j)} \in E(\overrightarrow{G})$ then $\overrightarrow{(v_i,v_j)}$ can be labeled (uniquely) as $k = (i-1)n-\binom{i}{2}+j-i \in J_N$. Moreover, if $k \in J_N$ is a label of an edge $\overrightarrow{e} \in E(\overrightarrow{G})$ then $\overrightarrow{e} = \overrightarrow{(v_i,v_j)}$.
\end{corollary}

Note that, $K_n$ is an underlying graph of $\overrightarrow{K_n}$. Therefore if $\overrightarrow{G}$ is a (directed) subgraph of $\overrightarrow{K_n}$ then its underlying graph $G$ is a subgraph of $K_n$.
Thus using Corollary \ref{t2}, for a given finite simple graph $G$ on $n$ unisolated vertices, we get an edge labeling for the graph $G$, which is uniquely determined.

\begin{theorem} \label{t3}
A digraph $\overrightarrow{K_n}$ is antimagic for $n \geq 3$.
\end{theorem}
\begin{proof}
Let $V(\overrightarrow{K_n}) = \{v_1, v_2, \ldots , v_n\}$ and let $E(\overrightarrow{K_n}) = \{\overrightarrow{e_k} = 
\overrightarrow{(v_i, v_j)} | k = (i-1)n - \binom{i}{2} + j - i, 1 \leq i < j \leq n\}$.
Let $f : E(\overrightarrow{K_n}) \to J_N$ be defined as $f(\overrightarrow{e_k}) = k$. Then by Theorem \ref{t1}, $f$ is bijective.

Now $S_{v_i}^- = \displaystyle \sum_{\overrightarrow{e} \in N^-(v_i)} f(\overrightarrow{e})$, where 
$N^-(v_i) = \{\overrightarrow{(v_p,v_i)} \in E(\overrightarrow{K_n}) | p < i \}$. 
Therefore $S_{v_i}^- = \displaystyle \sum_{p=1}^{i-1} f(\overrightarrow{(v_p,v_i)}) = \sum_{p=1}^{i-1} ((p-1)n - \binom{p}{2} + i - p)
= (n+1) \binom{i}{2} - n(i-1) - \frac{i(i-1)(i-2)}{6}$.

Similarly, $S_{v_i}^+ = \displaystyle \sum_{\overrightarrow{e} \in N^+(v_i)} f(\overrightarrow{e})$, where 
$N^+(v_i) = \{\overrightarrow{(v_i,v_q)} \in E(\overrightarrow{K_n}) | i < q \}$.
Therefore $S_{v_i}^+ = \displaystyle \sum_{q=i+1}^{n} f(\overrightarrow{(v_i,v_q)}) = \sum_{q=i+1}^{n} ((i-1)n - \binom{i}{2} + q - i)
= \binom{n}{2} + (n-(i+1)) (n(i-1) - \binom{i}{2})$.

But then the vertex sum of vertex $v_i$ is given by $S_{v_i} = S_{v_i}^- + S_{v_i}^+ = 
\frac{i^3}{3} - (n-1)i^2 + (n^2-n-\frac{4}{3})i - \frac{n(n-3)}{2}$.
We now claim that $S_{v_i} \neq S_{v_j}$ for $i \neq j$, i.e., $S_{v_i} \neq S_{v_j}$ for $v_i \neq v_j$.

For if, suppose $S_{v_i} = S_{v_j}$. Without loss of generality, suppose $i<j$.
As $S_{v_j} - S_{v_i} = 0$, we have $\frac{(j^3-i^3)}{3} - (n-1)(j^2-i^2) + (n^2-n-\frac{4}{3})(j-i) = 0$.
That is, $\frac{(j^2+ij+i^2)}{3} - (n-1)(j+i) + (n^2-n-\frac{4}{3}) = 0$, since $j-i>0$.
This implies that $n^2 - nx + y = 0$, where $x = i+j+1$ and $y = \frac{(i^2+ij+j^2-4)}{3} + (i+j)$.
Now $x^2 - 4y = \frac{19}{3} - 2(i+j) - \frac{(i-j)^2}{3} = 0$, if $i=1, j=2$, i.e., if $n = \frac{x}{2} = 2$, and $x^2 - 4y < 0$, otherwise.
This is a contradiction, since $n$ is an integer $\geq 3$.
Thus the vertex sums are pairwise distinct, and hence the proof.
\end{proof}

Note that in the proof of Theorem \ref{t3}, we get $x^2 - 4y < 0$ whenever $S_{v_i} = S_{v_j}$ for $i < j$ and $n \geq 3$. This implies that $n^2 - nx + y > 0$, and hence 
$S_{v_j} - S_{v_i} > 0$, i.e., $S_{v_i} < S_{v_j}$ whenever $i < j$.
Also $\overrightarrow{K_2}$ is not antimagic, and hence $K_2$ is not antimagic. 
Now $I(v_i) = N^-(v_i) \cup N^+(v_i)$ and $K_n$ is an underlying graph of $\overrightarrow{K_n}$. Therefore the following result immediately follows from Theorem \ref{t3}.

\begin{corollary}
A complete graph $K_n$ is antimagic for $n \geq 3$.
\end{corollary}

\begin{proposition}
Let $V(\overrightarrow{K_n}) = \{v_1, v_2, \ldots , v_n\}$ where $n \geq 3$. Then $S^-_{v_i} < S^-_{v_j}$ whenever $i < j$.
\end{proposition}
\begin{proof}
In the proof of Theorem \ref{t3}, we have seen that $S^-_{v_i} = (n+1) \binom{i}{2} - n(i-1) - \frac{i(i-1)(i-2)}{6}$. Therefore
$S^-_{v_j} - S^-_{v_i} = (n+1) (\binom{j}{2} - \binom{i}{2}) - n(j-i) - (\frac{j(j-1)(j-2)}{6} - \frac{i(i-1)(i-2)}{6})
= (n+1) (\frac{j^2-i^2}{2} - \frac{j-i}{2}) - n(j-i) - (\frac{(j^3-i^3)-3(j^2-i^2)+2(j-i)}{6})
= (j-i)((n+1)(\frac{i+j-1}{2}) - n - (\frac{(i^2+ij+j^2)-3(i+j)+2}{6}))
= (j-i)(n(\frac{i+j-3}{2}) + i+j - \frac{i^2+ij+j^2+5}{6})
= (j-i)(n(\frac{i+j-3}{2}) - \frac{(i+j-3)^2}{6} + \frac{ij+4}{6}) > 0$, since $1 \leq i < j$.
\end{proof}

\section{Total labelings of graphs/directed graphs} \label{sec3}

A {\it total labeling} of a graph $G$ is a bijection $f : V \cup E \to \{1, 2, \ldots , n+l\}$. A total labeling is called {\it super} if $f(V) = \{1, 2, \ldots , n\}$.
For a total labeling $f$, the associated {\it edge-weight} of an edge $e = \{u, v\}$ is defined by $w(f(e))= f(u)+f(v)+f(e)$.
The associated {\it vertex-weight} of a vertex $v$ is defined by $w_f(v) = f(v) + S_v$. 

We extend the above definitions to directed graphs as follows.
A {\it total labeling} of a directed graph $\overrightarrow{D}$ with vertex set $V$ and arc set $E$, is a bijection $f : V \cup E \to \{1, 2, \ldots , n+l\}$. A total labeling is called {\it super} if $f(V) = \{1, 2, \ldots , n\}$.
For a total labeling $f$, the associated {\it edge-weight} of a directed edge $\overrightarrow{e} = (\overrightarrow{u,v})$ is defined by $w(f(\overrightarrow{e}))= f(u)+f(v)+f(\overrightarrow{e})$.
The associated {\it vertex-weight} of a vertex $v$ is defined by $w_f(v) = f(v) + S^-_v + S^+_v$. 

A total labeling $f$ is called {\it edge-antimagic total (vertex-antimagic total)} if all edge-weights (vertex-weights) are pairwise distinct.
A graph (digraph) that admits an edge-antimagic total (vertex-antimagic total) labeling is called an {\it edge-antimagic total (vertex-antimagic total) graph (digraph)}.
A labeling that is simultaneously edge-antimagic total and vertex-antimagic total is called {\it totally antimagic total labeling}. 
A graph (digraph) that admits a totally antimagic total labeling is called a {\it totally antimagic total graph (digraph)}.
An edge-antimagic total labeling is called {\it super} if $f(V) = \{1, 2, \ldots , n\}$ and $f(E) = \{n+1, n+2, \ldots , n+l\}$.
A graph (digraph) that admits a super edge-antimagic total labeling is called a {\it super edge-antimagic total graph (digraph)}.

\begin{theorem} \label{t4}
For $n \geq 3$, $\overrightarrow{K_n}$ is vertex-antimagic total digraph.
\end{theorem}
\begin{proof}
Let $V(\overrightarrow{K_n}) = \{v_1, v_2, \ldots , v_n\}$ and let $E(\overrightarrow{K_n}) = \{\overrightarrow{e_k} = 
\overrightarrow{(v_i, v_j)} | k = (i-1)n - \binom{i}{2} + j - i, 1 \leq i < j \leq n\}$.
Let $f : V(\overrightarrow{K_n}) \cup E(\overrightarrow{K_n}) \to \{1, 2, \ldots , n+l\}$ be defined as $f(v_i) = i, \forall v_i \in V(\overrightarrow{K_n})$
and $f(\overrightarrow{e_k}) = n+k, \forall \overrightarrow{e_k} \in E(\overrightarrow{K_n})$. 
Then, by Theorem \ref{t1}, $f$ is bijective. Hence $f$ is not only total labeling but super also.
Therefore $w_f(v_i) = f(v_i) + S_{v_i} = i + \frac{i^3}{3} - (n-1)i^2 + (n^2-n-\frac{4}{3})i - \frac{n(n-3)}{2}
= \frac{i^3}{3} - (n-1)i^2 + (n^2-n-\frac{1}{3})i - \frac{n(n-3)}{2}$.

Suppose $v_i, v_j \in V(\overrightarrow{K_n})$ with $v_i \neq v_j$, i.e., $i \neq j$. Without loss, suppose $i < j$.
Consider $w_f(v_j) - w_f(v_i) = 0$. That is, $\frac{(j^3-i^3)}{3} - (n-1)(j^2-i^2) + (n^2-n-\frac{1}{3})(j-i) = 0$.
Therefore $\frac{(j^2+ij+i^2)}{3} - (n-1)(j+i) + (n^2-n-\frac{1}{3}) = 0$, since $j-i>0$.
This implies that $n^2 - nx + y = 0$, where $x = i+j+1$ and $y = \frac{(i^2+ij+j^2-1)}{3} + (i+j)$.
Now $x^2 - 4y = \frac{7}{3} - 2(i+j) - \frac{(i-j)^2}{3} = -(\frac{(i-j)^2 - 7}{3} + 2(i+j)) < 0$, since $1 \leq i < j \leq n$.
Thus $n^2 - nx + y > 0$, and hence $(j-i)(n^2 - nx + y) > 0$. Therefore $w_f(v_j) - w_f(v_i) > 0$, i.e., $w_f(v_i) < w_f(v_j)$.
\end{proof}

\begin{theorem} \label{t5}
Let $n \geq 3$. Let $\overrightarrow{e_k}, \overrightarrow{e_{k'}} \in E(\overrightarrow{K_n})$ with $\overrightarrow{e_k} \neq \overrightarrow{e_{k'}}$, i.e., $k \neq k'$.
If $k = (i-1)n - \binom{i}{2} + j - i$, $k' = (i'-1)n - \binom{i'}{2} + j' - i'$ with $n \neq \frac{2(j-j')}{(i'-i)} + \frac{i'+i-1}{2}$ whenever $1 \leq i<i'<j'<j \leq n$, then $\overrightarrow{K_n}$ is super edge-antimagic total digraph.
\end{theorem}

\begin{proof}
Let $V(\overrightarrow{K_n}) = \{v_1, v_2, \ldots , v_n\}$ and let $E(\overrightarrow{K_n}) = \{\overrightarrow{e_k} = 
\overrightarrow{(v_i, v_j)} | k = (i-1)n - \binom{i}{2} + j - i, 1 \leq i < j \leq n\}$.
Let $f : V(\overrightarrow{K_n}) \cup E(\overrightarrow{K_n}) \to \{1, 2, \ldots , n+l\}$ be defined as $f(v_i) = i, \forall v_i \in V(\overrightarrow{K_n})$
and $f(\overrightarrow{e_k}) = n+k, \forall \overrightarrow{e_k} \in E(\overrightarrow{K_n})$. 
Then by Theorem \ref{t1}, $f$ is bijective. Hence $f$ is not only total labeling but super also.
Therefore $w(f(\overrightarrow{e_k})) = f(v_i) + f(v_j) + f(\overrightarrow{e_k}) = i+j+(n+k) = i+j+n+((i-1)n - \binom{i}{2} + j - i) = ni + 2j - \binom{i}{2}$.

Suppose $\overrightarrow{e_k}, \overrightarrow{e_{k'}} \in E(\overrightarrow{K_n})$ with $\overrightarrow{e_k} \neq \overrightarrow{e_{k'}}$, i.e., $k \neq k'$. 
Without loss, suppose $k < k'$. Using Theorem \ref{t1}, suppose $k=F(i,j)$ and $k'=F(i',j')$ where $(i,j), (i',j') \in S$. 
Consider $w(f(\overrightarrow{e_{k'}})) - w(f(\overrightarrow{e_k})) = (ni' + 2j' - \binom{i'}{2}) - (ni + 2j - \binom{i}{2})
= n(i'-i) + 2(j'-j) - (\binom{i'}{2} - \binom{i}{2}) = 2(j'-j) + (i'-i)(n - \frac{i'+i-1}{2})$.
As $k < k'$, $(i,j) < (i',j')$. Therefore either $i < i'$, or $i=i'$ and $j < j'$.
If $i < i'$ then we have either $j \leq j'$ or $j > j'$.
Now if $i < i'$ with $j \leq j'$, and $i=i'$ with $j < j'$ then $w(f(\overrightarrow{e_{k'}})) - w(f(\overrightarrow{e_k})) > 0$, since 
$n - \frac{i'+i-1}{2} = \frac{n-i'}{2} + \frac{n-i}{2} + \frac{1}{2}$ and $i < i' \leq n$.
Also if $i < i'$ with $j > j'$ then we have $1 \leq i<i'<j'<j \leq n$ with $n \geq 4$.
Therefore $w(f(\overrightarrow{e_{k'}})) - w(f(\overrightarrow{e_k})) = 0$ if and only if $(i'-i)(n - \frac{i'+i-1}{2}) = 2(j-j')$.
Thus if $n \neq \frac{2(j-j')}{(i'-i)} + \frac{i'+i-1}{2}$, where $1 \leq i<i'<j'<j \leq n$, then $w(f(\overrightarrow{e_{k'}})) - w(f(\overrightarrow{e_k})) \neq 0$. 
Hence the proof.
\end{proof}

The following result follows from Theorem \ref{t4} and \ref{t5}.

\begin{corollary} \label{c4}
Let $n \geq 3$. Let $\overrightarrow{e_k}, \overrightarrow{e_{k'}} \in E(\overrightarrow{K_n})$ with $\overrightarrow{e_k} \neq \overrightarrow{e_{k'}}$, i.e., $k \neq k'$.
If $k = (i-1)n - \binom{i}{2} + j - i$, $k' = (i'-1)n - \binom{i'}{2} + j' - i'$ with $n \neq \frac{2(j-j')}{(i'-i)} + \frac{i'+i-1}{2}$ whenever $1 \leq i<i'<j'<j \leq n$, then $\overrightarrow{K_n}$ is totally antimagic total digraph.
\end{corollary}

\begin{corollary}
Let $n \geq 3$. Let $e_k, e_{k'} \in E(K_n)$ with $e_k \neq e_{k'}$, i.e., $k \neq k'$.
If $k = (i-1)n - \binom{i}{2} + j - i$, $k' = (i'-1)n - \binom{i'}{2} + j' - i'$ with $n \neq \frac{2(j-j')}{(i'-i)} + \frac{i'+i-1}{2}$ whenever $1 \leq i<i'<j'<j \leq n$, then $K_n$ is totally antimagic total graph.
\end{corollary}

\section{Antimagic orientation of graphs} \label{sec4}

Recall that, an {\it antimagic labeling of a directed graph} $D$ with $n$ vertices and $l$ arcs (or directed edges) is a bijection from the set of arcs of $D$ to the integers $\{1, 2, \ldots , l\}$ such that all $n$ oriented vertex sums are pairwise distinct, where an {\it oriented vertex sum} is the sum of labels of all edges entering that vertex minus the sum of labels of all edges leaving it. 
Hefetz et al.\ \cite{HMS} raised two questions: Is every orientation of any simple connected undirected graph antimagic?,  
Given any undirected graph $G$, does there exist an orientation of $G$ which is antimagic? They call such an orientation an {\it antimagic orientation} of $G$.
In the following, we prove that there exists an antimagic orientation of $K_n$.
\begin{theorem}
For $n \geq 3$, $K_n$ has an antimagic orientation.
\end{theorem}
\begin{proof}
Let $V(\overrightarrow{K_n}) = \{v_1, v_2, \ldots , v_n\}$ and let $E(\overrightarrow{K_n}) = \{\overrightarrow{e_k} = 
\overrightarrow{(v_i, v_j)} | k = (i-1)n - \binom{i}{2} + j - i, 1 \leq i < j \leq n\}$.
Let $f : E(\overrightarrow{K_n}) \to J_N$ be defined as $f(\overrightarrow{e_k}) = k$. Then by Theorem \ref{t1}, $f$ is bijective.

Now $S_{v_i}^- = \displaystyle \sum_{\overrightarrow{e} \in N^-(v_i)} f(\overrightarrow{e})$, where 
$N^-(v_i) = \{\overrightarrow{(v_p,v_i)} \in E(\overrightarrow{K_n}) | p < i \}$. 
Therefore $S_{v_i}^- = \displaystyle \sum_{p=1}^{i-1} f(\overrightarrow{(v_p,v_i)}) = \sum_{p=1}^{i-1} ((p-1)n - \binom{p}{2} + i - p)
= (n+1) \binom{i}{2} - n(i-1) - \frac{i(i-1)(i-2)}{6}$.

Similarly, $S_{v_i}^+ = \displaystyle \sum_{\overrightarrow{e} \in N^+(v_i)} f(\overrightarrow{e})$, where 
$N^+(v_i) = \{\overrightarrow{(v_i,v_q)} \in E(\overrightarrow{K_n}) | i < q \}$.
Therefore $S_{v_i}^+ = \displaystyle \sum_{q=i+1}^{n} f(\overrightarrow{(v_i,v_q)}) = \sum_{q=i+1}^{n} ((i-1)n - \binom{i}{2} + q - i)
= \binom{n}{2} + (n-(i+1)) (n(i-1) - \binom{i}{2})$.

But then the oriented vertex sum of vertex $v_i$ is given by $S^o_{v_i} = S_{v_i}^- - S_{v_i}^+ = -\frac{2}{3} i^3 + (2n+1)i^2 - (n^2+2n+\frac{1}{3})i + \binom{n+1}{2}$.
We claim that $S^o_{v_i} \neq S^o_{v_j}$ for $v_i \neq v_j$, i.e., $S^o_{v_i} \neq S^o_{v_j}$ for $i \neq j$. 

Without loss, suppose $i<j$. Now suppose $S^o_{v_j} - S^o_{v_i} = 0$. That is, $-\frac{2}{3} (j^3-i^3) + (2n+1)(j^2-i^2) - (n^2+2n+\frac{1}{3})(j-i) = 0$. Therefore we have
$\frac{2}{3} (i^2+ij+j^2) - (2n+1)(i+j) + (n^2+2n+\frac{1}{3}) = 0$, since $j-i \neq 0$.
That is, $n^2-nx+y = 0$, where $x = 2(i+j-1)$ and $y = \frac{1 - 3(i+j) + 2(i^2+ij+j^2)}{3}$. Now $x^2-4y = \frac{4}{3} ((i+j)^2 + i(j-1) + j(i-1) + 2) > 0$. 
Hence $n = (i+j-1) \pm \sqrt{\frac{1}{3} ((i+j)^2 + i(j-1) + j(i-1) + 2)}$. 
Thus $n$ is dependent on $i$ and $j$. That is, $n$ is not unique, a contradiction. Therefore $S^o_{v_i} \neq S^o_{v_j}$. Hence the proof.
\end{proof}

\end{document}